\documentclass[10pt]{article}
\usepackage{lineno,hyperref}

\usepackage{geometry}
\usepackage{graphicx}
\usepackage{amssymb}
\usepackage{amsmath}
\usepackage{amsthm}
\usepackage{enumerate}
\usepackage{here}
\usepackage{setspace}

  \newtheorem{Theorem}{Theorem}[section]
  \newtheorem{Lemma}[Theorem]{Lemma}

  \newtheorem{Property}[Theorem]{Property}
  \newtheorem{Definition}[Theorem]{Definition}
  \newtheorem{Example}[Theorem]{Example}
  
\newcommand{\Proof}{\noindent{\bf Proof}\quad}
\newcounter{nombre}

\title{\textbf{GDD type Spanning Bipartite Block Designs}}
\author{Shoko Chisaki
\thanks{Department of Information Systems, Osaka Institute of Technology}
\and Ryoh Fuji-Hara
\thanks{Faculty of Engineering, Information and Systems, University of Tsukuba}
\and Nobuko Miyamoto
\thanks{Department of Information Sciences, Tokyo University of Science}}
\date{}

\begin{document}

\maketitle

\begin{abstract}
There is a one-to-one correspondence between the point set of a group divisible design (GDD) with $v_1$ groups of $v_2$ points and the edge set of a complete bipartite graph $K_{v_1,v_2}$.
A block of GDD corresponds to a subgraph of $K_{v_1,v_2}$.
A set of subgraphs of  $K_{v_1,v_2}$ is constructed from a block set of GDDs.
If the GDD satisfies the  $\lambda_1, \lambda_2$  concurrence  condition, then the set of subgraphs
 also satisfies the spanning bipartite block design (SBBD)  conditions \cite{Chisa3}.
We also propose a method to construct SBBD directly from an $(r,\lambda)$-design and a difference matrix over a group. Suppose the $(r,\lambda)$-design consists of  $v_2$ points and $v_1$ blocks. When $v_1 >> v_2$, we show a method to construct a SBBD with $v_1$ is close to $v_2$ by partitioning the block set.
\end{abstract}

\section{Introduction}

Let $V_1=\{1,2,\ldots,v_1\}$ and $V_2=\{1,2,\ldots,v_2\}$ be disjoint two  point sets, and $E=\{ e_{ij} \, | \, i\in V_1, j\in V_2\}$ be the edge set between $V_1$ and $V_2$. $K_{v_1,v_2} = (V_1 , V_2\, ; \, E)$ is the complete bipartite graph with $v_1, v_2$ point sets.
Let $\mathcal{B}=\{B_1,B_2,\ldots, B_N \}$ be a collection of subgraphs of $K_{v_1,v_2}$  called spanning bipartite blocks (SB-blocks).
If $\mathcal{B}$ satisfies the following five conditions, then we call $(K_{v_1,v_2}\ ;\  \mathcal{B} )$ a \textit{spanning bipartite block design} (SBBD):

\begin{enumerate}
\item   Each SB-block $B_i$ of $\mathcal{B}$ is incident with all points of $V_1$ and  $V_2$ (spanning condition).
\item  Each edge of $K_{v_1,v_2}$ appears in $\mathcal{B}$ exactly $\mu$ times.
\item   Any  two edges $e_{ij}, e_{ij'}$ such that $i \in V_1$, $j, j' \in V_2, (j \ne j')$ are contained together in $\lambda_{12}$ SB-blocks in $\mathcal{B}$.
\item     Any two edges $e_{ij}, e_{i'j}$ such that  $i, i' \in V_1, (i \ne i')$, $j \in V_2$ are contained  together  in $\lambda_{21}$ SB-blocks in $\mathcal{B}$.
\item     Any  two edges $e_{ij}$, $e_{i'j'}$ such that $i, i' \in V_1, (i \ne i')$, \, $j, j'\in V_2,  (j \ne j')$ are contained together in $\lambda_{22}$ SB-blocks in $\mathcal{B}$.
\end{enumerate}

A spanning bipartite block design is first proposed in Chisaki et al.~\cite{Chisa3}.
This design is for a statistical model to estimate treatment parameters with the structure of a complete bipartite graph.
In \cite{Chisa3}, it is proved that SBBDs with some conditions is A-optimum.
SBBD can also be used as a kind of sparsing method to prevent over-fitting  of deep learning.
Compared to the random dropconnect method in Wan et al.~\cite{wan2013}, the spanning condition plays an important role in sparsing neural networks  to drop out connections independently at each layer.
And the balancing properties work for reducing the variances of weight estimators.

There is a similar block design called
a {\it balanced bipartite block design} (BBBD).
A BBBD $(V_1,V_2 : \mathcal{B})$ is defined belows:
\begin{enumerate}
\item[1:] let $K_{v_1,v_2}$ be a complete bipartite graph with $V_1$, and $V_2$ point sets, $|V_1|=v_1,|V_2|=v_2$,
\item[2:]  $\mathcal{B} =\{ B_1, B_2,\ldots, B_N \}$ be a set of complete sub-bipartite graphs $K_{k_1,k_2}$ of $K_{v_1,v_2}$ (called blocks) and block size is $k=k_1+k_2$.
\item[3:] for any $t_1$ points from $V_1$ and $t_2$ points from $V_2$, there are exactly $\mu_{t_1,t_2}$ blocks in $\mathcal{B}$ containing those points.
\end{enumerate}
Although this design is similar in name to SBBD, even the  blocks of BBBD are all complete sub-bipartite graphs $K_{k_1,k_2}$.
In Chisaki et al.~\cite{Chisa1, Chisa2}, different designs from SBBD were proposed for the same purpose of deep learning. Those designs are rather close to BBBD.
Many papers, including
Kageyama and Sinha~\cite{Kageyama-Sinha1988}, Mishima et al.~\cite{Mishima2001} and Jaggi et al.~\cite{Jaggi-etal-2023}  show constructions  satisfying the third condition for $\mu_{2,0}, \mu_{0,2}$ and $\mu_{1,1}$.
Ozawa et al.~\cite{ozawa2002} show  constructions of BBBD (it is  called a {\it split-block designs} in \cite{ozawa2002} )
satisfying the third condition of $\mu_{t_1, t_2}$ for $0 \le t_1, t_2 \le 2$.
Martin~\cite{Martin1998} defined a type of BBBD (it is called {\it Mixed $t$-design})
satisfying the third condition for any $t_1,t_2$ such that $t_1+t_2=t$.
He shows some constructions for $t=2$ and $t=3$.
\\[0.1cm]
In this paper, we show that an SBBD with a certain condition and a GDD can be considered equivalent,
and propose a construction for the SBBDs using an $(r,\lambda)$-design and a difference matrix.
Additionally, we describe the E-optimality of SBBDs and show some examples.

\section{Design matrix}
In this section, we introduce a matrix representation of SBBDs and re-express the five conditions.
First, we define a $(0,1)$-matrix $X$ from the SB-blocks called a {\it design matrix}.

\begin{itemize}
\item Suppose that the edges $e_{ij}$ of $K_{v_1,v_2}$ are arranged in the following lexicographical order:
\begin{equation*}
( e_{11}, e_{12},\ldots, e_{1v_2}\ ;\  e_{21},e_{22},\ldots, e_{2v_2}\ ; \   \cdots \ ;\ e_{v_11},\ldots, e_{v_1v_2}).
\end{equation*}
This sequence of edges corresponds to the columns of $X$.
Denote $(e_{ij})$ for the column number corresponding to the edge $e_{ij}$.
\item Put $X=[x_{k,(e_{ij})}]$, then  $x_{k,(e_{ij})}$  is the element of the $k$-th row and the $(e_{ij})$-th column of $X$.
The design matrix $X$ is defined by the SB-blocks  $B_1,B_2,\ldots,B_N$ as follows:
$$x_{k,(e_{ij})}=
\begin{cases}
1 &\mbox{ if}\ \   e_{ij} \in B_k \\
0 &\mbox{ otherwise}
\end{cases}
$$
\item $X$ is an $(N \times v_1v_2)$-matrix.
\end{itemize}

Let $X_i$ be an $(N \times v_2)$-submatrix consisting of $v_2$ columns of $X$ corresponding to $(e_{i1}, e_{i2},\ldots,$ $e_{iv_2})$.
Then the design matrix $X$ is partitioned into $v_1$ submatrices expressed as $X=( X_1 | X_2 | \cdots$  $|X_{v_1} )$.
If $(K_{v_1,v_2}\,;\, \mathcal{B} )$ is a spanning bipartite block design then
$X=( X_1 | X_2 | \cdots | X_{v_1})$ has the following property:
\begin{enumerate}[(i)]
\item any row of $X_i$ is not  zero-vector for $1 \leq i \leq v_1$ and
 $\sum_{i=1}^{v_1} X_i$  does not contain a zero element (spanning condition),
\item diag$(X_i ^t X_i) = (\mu, \mu,\ldots, \mu)$ for $1\le i\le v_1$,
\item all off-diagonal elements of $X_i^{\,t} X_i$ are $\lambda_{12}$ for $1 \leq i \leq v_1$,
\item  diag$(X_i^{\,t} X_j)=(\lambda_{21},\lambda_{21}, \ldots, \lambda_{21}) $ for  $1 \leq i \ne j \leq v_1$,
\item all off-diagonal elements of $X_i^{\,t} X_j$ are $\lambda_{22}$ for $1 \leq i \ne j \leq v_1$.
\end{enumerate}

$X^tX$ is called an \textit{information matrix}.
The information matrix of SBBD is expressed as follows:
\begin{align*}
X^t X &=I_{v_1}\otimes (X_i^t X_i) +(J_{v_1}-I_{v_1})\otimes (X_i^t X_j ) \\
&=
 I_{v_1} \otimes
\left[
\begin{array}{cccc}
\mu & \lambda_{12} & \cdots & \lambda_{12}\\
\lambda_{12} & \mu & \cdots & \lambda_{12}\\
\vdots &  \vdots &   \ddots           &  \vdots \\
\lambda_{12} & \lambda_{12} & \cdots & \mu
\end{array}
\right]
+ (J_{v_1}-I_{v_1}) \otimes
\left[
\begin{array}{cccc}
\lambda_{21} & \lambda_{22} & \cdots & \lambda_{22}\\
\lambda_{22} &\lambda_{21} & \cdots & \lambda_{22}\\
\vdots &  \vdots &   \ddots          &  \vdots \\
\lambda_{22} & \lambda_{22} & \cdots & \lambda_{21}
\end{array}
 \right],
\end{align*}
where $I_n$ is the identity matrix of size $n$ and $J_n$ is the $(n \times n)$  all-ones matrix.
A matrix expressed by $a I_n + b (J_n-I_n)$ is called {\it completely symmetric}.
The information matrix above has a double structure of a completely symmetric matrix.
The spanning bipartite block design is denoted as SBBD($v_1,v_2, N ; \Lambda$), where $\Lambda= (\mu, \lambda_{12}, \lambda_{21}, \lambda_{22})$.

\begin{Example}\label{ex2}
Let
\begin{center}
$X = (X_1 | X_2 | X_3 ) =$
\scalebox{0.8}[0.8]{
$
\left[
\begin{array}{ccc|ccc|ccc}
 0 & 1 & 1 & 1 & 1 & 0 & 1 & 1 & 0 \\
 1 & 0 & 1 & 0 & 1 & 1 & 0 & 1 & 1 \\
 1 & 1 & 0 & 1 & 0 & 1 & 1 & 0 & 1 \\
 0 & 1 & 1 & 0 & 1 & 1 & 1 & 0 & 1 \\
 1 & 0 & 1 & 1 & 0 & 1 & 1 & 1 & 0 \\
 1 & 1 & 0 & 1 & 1 & 0 & 0 & 1 & 1 \\
 0 & 1 & 1 & 1 & 0 & 1 & 0 & 1 & 1 \\
 1 & 0 & 1 & 1 & 1 & 0 & 1 & 0 & 1 \\
 1 & 1 & 0 & 0 & 1 & 1 & 1 & 1 & 0 \\
\end{array}
\right]
$
}
\end{center}
be a design matrix of an SBBD.
Then the information matrix is

$$ X^t X=
I_3 \otimes
\left[\begin{array}{ccc}6 & 3 & 3 \\3 & 6 & 3 \\3 & 3 & 6\end{array}\right]
+ (J_3-I_3)\otimes  \left[\begin{array}{ccc}4 & 4 & 4 \\4 & 4 & 4 \\4 & 4 & 4\end{array}\right].
$$
 The design matrix $X$ satisfies the spanning condition since any  row  of  $X_i$ is no zero-vector, and $X_1+X_2+X_3$ does not contain $0$.
So we have an SBBD($3, 3, 9 ; \Lambda),$  $\Lambda=(6,3,4,4)$.
\end{Example}

As you can see from the above example, the spanning condition can not be confirmed from the information matrix $X^t X$.
If  $v_1 \ll v_2$, there is a high possibility that the spanning condition is not met.
Such a design in which the spanning condition is not guaranteed is denoted by SBBD$^*$.

\section{Group Divisible Designs and SBBDs}

\begin{Definition}[Group Divisible Design, see Beth et al.~\cite{Beth1999}]\label{def:gdd}
Let $V$ be the $v$-point set which is
partitioned into $G_1,G_2,\ldots, G_m$, called {\it groups},
and $\mathcal{B}=\{B_1,B_2,\ldots, B_N\}$ ({\it blocks}) is a collection of subsets of $V$.
If $(V , \mathcal{B})$ satisfies the following conditions, it is called a {\it group divisible design} or simply GDD:
\begin{enumerate}
\item[(1)] any pair of distinct two points in the same group is contained in precisely $\lambda_1$ blocks.
\item[(2)] any pair of two points in distinct groups is contained in precisely $\lambda_2$ blocks.
\end{enumerate}
In this paper, we add the following two conditions:
\begin{enumerate}
\item[(3)] each group has the same number of points, $|G_i |= g$, for $i=1,2,\ldots, m$, i.e. $v=m g $,
\item[(4)] each point of $V$ is contained in exactly  $r$ blocks, i.e. $r = (\sum_{i=1}^N | B_i |) / v$
\end{enumerate}
It is denoted by GD$_{\lambda_1,\lambda_2}(K , g \,;\, v)$, where $K$ is the set of block sizes, or by GD$_{\lambda_1,\lambda_2}(k , g\, ; \,  v)$ if $K=\{k\}$.
A GD$_{0,\lambda}(m, g\, ;\,m g)$ is said to be a {\it transversal design} or an  {\it orthogonal array}.
\end{Definition}

\begin{Property}[Bose and Connor~\cite{Bose1952}]
The parameters of GD$_{\lambda_1,\lambda_2} (k,  g \,;\, v)$ with $N$ blocks and $v=mg$ have the following relation:
\begin{equation*}
 k N  = v r ,\ \  (g-1) \lambda_1+g(m-1)\lambda_2 = r (k-1),\ \  r \ge \lambda_1, \lambda_2.
 \end{equation*}
 \end{Property}

Let $V_1=\{1,2,\ldots,v_1\}$ and $V_2 =\{1,2,\ldots, v_2\}$.
Consider the complete bipartite graph $K_{v_1,v_2}=(V_1, V_2 ; E)$,
where the edge set is $E=\{ e_{ij}\, |\ i\in V_1, j\in V_2\}$.
Let $V=\{p_{11}, p_{12},\ldots, p_{v_1v_2}\}$ be  the point set of a $GD_{\lambda_1,\lambda_2}(k, v_2; v_1v_2)$   with $v_1$ groups, where  $G_i=\{p_{i1},p_{i2},\ldots, p_{iv_2}\}$.
Then there is  a  one-to-one correspondence between the point set  $V$
and the edge set $E$ of $K_{v_1,v_2}$ such as:
\begin{equation*}
p_{ij}\in V  \Leftrightarrow  e_{ij} \in E .
\end{equation*}
From this correspondence, a block of GDD is considered an SB-block.
A GDD satisfies the conditions of SBBD except for the spanning condition.
We can easily  see the following result:
\begin{Property}
If $(V, \mathcal{B} )$ is a GD$_{\lambda_1,\lambda_2}(K,v_2\,;\,v_1v_2)$, then it is also an SBBD*$(v_1,v_2,N ; \Lambda)$, $\Lambda=(\mu,\lambda_{12},\lambda_{21},$ $\lambda_{22})$ with the following relations:
\begin{equation*}
 r = \mu \ , \ \lambda_1 =\lambda_{12}\ , \ \lambda_2 =\lambda_{21} = \lambda_{22}.
\end{equation*}
\end{Property}

If a GDD satisfies the following conditions,
the SBBD$^*$ is an SBBD:
\begin{itemize}
    \item For every block $B \in \mathcal{B}$ and every group $G_i$, $|B \cap G_i|\ge 1$,
    \item Every element of $V_2$ appears at least once in  the set of the second subscripts of points in $B$ for every block $B \in \mathcal{B}$, i.e. $\{ j \ |\ p_{ij}\in B\}= V_2$.
\end{itemize}

A GDD not satisfying the second condition may be able to adjust to satisfy the spanning condition using the following property:
\begin{Property}
Let $\delta$ be a permutation on $\{1,2,\ldots, v_2\}$.
Even if the points within a  group $G_i=\{p_{i1}, p_{i2},\ldots,$ $p_{i v_2} \}$ are  rearranged by $\delta$ as:
$$\{p_{i\delta(1)}, p_{i\delta(2)}, \ldots, p_{i\delta(v_2)}\},$$
they remain a GDD with the same parameters.
\end{Property}
An SBBD satisfying $\lambda_{21}=\lambda_{22}$ is called a {\it GDD-type} SBBD.

\begin{Example} \label{ex:3}
Consider $GD_{3,4}(6,3\, ; 9), N=9$.
The points of the groups are represented here as $G_1=\{ 1_1,1_2,1_3\},$  $\  G_2=\{ 2_1, 2_2, 2_3\},
\  G_3=\{ 3_1,3_2,3_3 \}$,  and the blocks are:
\begin{equation}
\begin{array}{ccc}
B_1= \{1_2, 1_3\ ;\  2_2, 2_3\ ;\ 3_2, 3_3\}, &
B_2= \{ 1_1, 1_3\ ;\  2_1, 2_3 \ ;\ 3_1, 3_3 \}, &
B_3= \{ 1_1, 1_2\ ;\  2_1, 2_2 \ ;\ 3_1, 3_2 \},
 \\
B_4=\{ 1_1, 1_3\ ;\ 2_1, 2_2 \ ;\ 3_2, 3_3 \}, &
B_5=\{ 1_2, 1_3\ ;\ 2_1, 2_3 \ ;\ 3_1, 3_2 \}, &
B_6=\{ 1_1, 1_2\ ;\ 2_2, 2_3 \ ;\ 3_1, 3_3 \},
 \\
B_7=  \{ 1_1, 1_2\ ;\ 2_1, 2_3 \ ;\ 3_2, 3_3 \}, &
B_8=  \{ 1_2, 1_3\ ;\ 2_1, 2_2 \ ;\ 3_1, 3_3 \}, &
B_9=  \{ 1_1, 1_3\ ;\ 2_2, 2_3 \ ;\ 3_1, 3_2 \}.
\end{array} \nonumber
\end{equation}
This is from AG($2,3$),
 the group is a  parallel class of lines, and the blocks are the complement  of lines that transverse the parallel lines.
Let $\psi(B)=\{j \,|\, i_j \in B\}$.
$\psi(B_1)$ is missing 1,  $\psi(B_2)$ is missing 2 and $\psi(B_3)$ is missing 3.
This does not satisfy the spanning conditions.
By a cyclic permutation $\delta=(123)$ on the subscripts of  $G_3$ points, i.e.
$ 3_1\mapsto 3_2,\ 3_2\mapsto 3_3,\  3_3\mapsto 3_1$,
we have the following GDD:
\begin{equation}
\begin{array}{ccc}
B_1= \{1_2, 1_3\ ;\  2_2, 2_3 \ ;\ 3_3, 3_1 \}, &
B_2= \{1_1, 1_3\ ;\  2_1, 2_3 \ ;\ 3_2, 3_1 \}, &
B_3= \{1_1, 1_2\ ;\  2_1, 2_2 \ ;\ 3_2, 3_3 \},
 \\
B_4=\{1_1, 1_3\ ;\ 2_1, 2_2 \ ;\ 3_3, 3_1 \}, &
B_5=\{1_2, 1_3\ ;\ 2_1, 2_3 \ ;\ 3_2, 3_3 \}, &
B_6=\{1_1, 1_2\ ;\ 2_2, 2_3 \ ;\ 3_2, 3_1 \},
 \\
B_7=  \{1_1, 1_2\ ;\ 2_1, 2_3 \ ;\ 3_3, 3_1 \}, &
B_8=  \{1_2, 1_3\ ;\ 2_1, 2_2 \ ;\ 3_2, 3_1 \}, &
B_9=  \{1_1, 1_3\ ;\ 2_2, 2_3 \ ;\ 3_2, 3_3 \}.
\end{array} \nonumber
\end{equation}
Their information matrices are both
\begin{equation*}
 \mathbf{X}^t \mathbf{X}=
I_3 \otimes
\left[\begin{array}{ccc}6 & 3 & 3 \\3 & 6 & 3 \\3 & 3 & 6\end{array}\right]
+ (J_3-I_3)\otimes  \left[\begin{array}{ccc}4 & 4 & 4 \\4 & 4 & 4 \\4 & 4 & 4\end{array}\right].
\end{equation*}
The second example  is a GDD-type  SBBD$(3, 3, 9 ; \Lambda),$  $\Lambda=(6,3,4,4)$.
\end{Example}

\section{Construction from an $(r,\lambda)$-design and a Difference Matrix}

In this section, we show a construction of GDD type SBBD that is not from a group divisible design using an $(r,\lambda)$-design and a difference matrix.
Our idea for constructing SBBD consists of the following three steps:
\begin{description}
\item[First:  ] We select an incidence matrix $H$ of an $(r,\lambda)$-design,
\item[Second: ] Using the incidence matrix $H$ as a seed, a set of  matrices called tile matrices are  generated by the operation of a group,
\item[Third:  ] A design  matrix $X$ of SBBD can be constructed by pasting the tile matrices on a combinatorial array called a difference matrix over the group.
\end{description}

\begin{Definition}[$(r,\lambda)$-design, Stanton and Mullin~\cite{StantonMullin1966}]
Let $V$ be a $v$-point set and $\mathcal{B}$ a collection of subsets (blocks) of $V$.
If $(V, \mathcal{B})$ holds the following conditions, it is called an $(r,\lambda)$-design with $v$ points and  $b$ blocks:
  \begin{itemize}
  \item each point of $V$ is contained in exactly $r$ blocks of $\mathcal{B}$,
  \item any two distinct points of $V$ are contained in precisely $\lambda$ blocks of $\mathcal{B}$.
  \end{itemize}
\end{Definition}

If the block size is a constant $k$ for each block, then it is called a \textit{balanced incomplete block design}, denoted  by $(v,k,\lambda)$-BIBD, and if $|V| = |\mathcal{B}|$, it is called \textit{symmetric design}, then $r=k$.

\begin{Definition}[Difference Matrix over a group $\mathbf{E}_b$, Jungnickel~\cite{Jungnickel1979}]
Le $D=[d_{ij}]$ be an $(\eta b \times s)$-matrix on a group $\mathbf{E}_b$ of order $b$, and
 $D(i,j)=\{(d_{ki},d_{kj})\,|\,k=1,2,\ldots, \eta v\}$, $1\le i\ne j \le s$.
If the multi-set
$$\{d-d' \, |\,  (d,d') \in D(i,j)\,  \}$$
contains each element of $\mathbf{E}_b$ precisely $\eta$ times for any $1 \le i\ne j \le s$, then
$D$ is called a $(b, s ; \eta)$-difference matrix (DM) over $\mathbf{E}_b$.
\end{Definition}
If $s=b\eta$, then $D$ may be  called a \textit{generalized Hadamard matrix}.
On difference matrices, we have the following well-known properties, see Beth et al.~\cite{Beth1999}:

\begin{Property}[Beth et al.~\cite{Beth1999}]\label{pr:add}
Let $D$ be a difference matrix.
A matrix $D'$ obtained by adding  an element $c \in \mathbf{E}_b$ to all elements of a column of $D$ is also a  difference matrix.
$$ D' = [d'_{ij}] \mbox{ such that } d'_{ij}\equiv d_{ij}+ c \mbox{ \  for } i=1,2,\ldots, \eta b.$$
\end{Property}

Using this property, it can be  adjusted to satisfy  the spanning condition of SBBD.

\begin{Property}[Beth et al.~\cite{Beth1999}]\label{ppMD}
For any prime power $q$, there exists a $(q, q ; 1)$-DM.
\end{Property}
Many examples of existence, such as $(r,\lambda)$-designs and difference matrices, are shown in Colbourn and Dinitz~\cite{Handbook}.

Let $(V,\mathcal{B})$ be an $(r,\lambda)$-design with $v$ points and $b$ blocks,

\begin{equation*}\label{eq:H}
H=
 \begin{bmatrix}
 \mathbf{h}_{x_0}\\
 \mathbf{h}_{x_1}\\
 \vdots\\
  \mathbf{h}_{x_{b-1}}
 \end{bmatrix},
\end{equation*}
where $\mathbf{h}_{x_0}, \mathbf{h}_{x_1},\ldots, \mathbf{h}_{x_{b-1}}$ are the row vectors of $H$  and their subscripts are described by the elements of $\mathbf{E}_b=\{x_0, x_1,\ldots, x_{b-1} \}$ arranged in a certain order.
The tile matrix $T_{y}$ is an array of rows  rearranged by adding the element $y$ of $\mathbf{E}_b$ to the subscripts of each row of $H$ as follows:

\begin{equation}\label{eq:T}
T_{y} =
\begin{bmatrix}
\mathbf{h}_{x_0 + y}\\
\mathbf{h}_{x_1 + y}\\
\vdots\\
\mathbf{h}_{x_{b-1} + y}
\end{bmatrix} \mbox{ for }  y \in \mathbf{E}_b.
\end{equation}

 Assume $x_0 = 0$ (identity) in  $\mathbf{E}_b$, that is, $T_{x_0} =H$.
Each $T_y$ has the following properties:
\begin{itemize}
\item $T_y$ is a $(b \times v)$-matrix for $y \in \mathbf{E}_b$,
\item the set of rows of  $T_y$  is precisely equal to the set of rows of $H$,
which implies
\begin{gather}
 T_y^t \, T_y = H^t H = rI_v +\lambda (J_v-I_v),  \mbox{  for any } y \in \mathbf{E}_b , \label{pr:TtT}
 \end{gather}
 \end{itemize}
 Then we have following equations about the tile matrix $T_y$:
\begin{Lemma}\label{lem:diff}
 For any $x, y , d \in \mathbf{E}_b$,
 it holds
 \begin{equation}\label{eq:diff}
     T_x^t \, T_y = (T_{x+d})^t \, T_{y+d}\,.
 \end{equation}
 For any $x \in \mathbf{E}_b$,
  it holds
\begin{equation}\label{eq:diff3}
    \sum_{y \in \mathbf{E}_b} T_x^t \, T_y = r^2 J_v \,.
\end{equation}
\end{Lemma}

\Proof
Let $\mathbf{E}_b = \{x_0, x_1,\ldots, x_{b-1}\}$ be a group of order $b$.
Let a pair of the  $x_i$-th rows of $T_x$ and $T_y$ be ($\mathbf{h}_{x_i+x}, \mathbf{h}_{x_i+y}$).
Similar pair of $T_{x+d}$ and $T_{y+d}$ is described as ($\mathbf{h}_{y_i+x+d}, \mathbf{h}_{y_i+y+d}$), $y_i \in \mathbf{E}_b$.
If $y_i= x_i-d$, then those two pairs are equal. That is, the set of pairs \{($\mathbf{h}_{x_i+x}, \mathbf{h}_{x_i+y})\, ; \, x_i \in \mathbf{E}_b\}$  is the same as the set of pairs  \{($\mathbf{h}_{y_i+x+d}, \mathbf{h}_{y_i+y+d})\, ; \, y_i \in \mathbf{E}_b\}$, which implies that $$T_x^t \, T_y = (T_{x+d})^t \, T_{y+d}\,.$$
Next, it is easy to see that any row of $\sum_{y \in \mathbf{E}_b} T_y$ equals to $\sum_{i=0}^{b-1} \mathbf{h}_{x_i}$, and therefore
$\sum_{y \in \mathbf{E}_b} T_y= r J_{b,v}$.
Hence we have $$\sum_{y \in \mathbf{E}_b} T_x^t \, T_y = r^2 J_v\, .$$
\qed

Let $D=[d_{i,j}]$ be an
$(\eta b \times s)$-matrix of $(b,s ; \eta)$-DM over $\mathbf{E}_{b}$.
We paste the tile matrices $T_0,T_{x_1},\ldots, T_{x_{b-1}}$  on $D$  to make a design matrix
\begin{equation}\label{def_X}
X = [T_{d_{i,j}}]=(X_1 | X_2 | \cdots | X_s).
\end{equation}

This $X$ is an $(\eta b^2 \times s v )$-matrix, and
\begin{equation}\label{def_Xj}
X_j =
\begin{bmatrix}
T_{d_{1,j}}\\
T_{d_{2,j}}\\
\vdots\\
T_{d_{\eta b, j}}\\
\end{bmatrix}\quad \textrm{for } 1\le j\le s.
\end{equation}
We have the next theorem regarding each row of $X$ as a new SB-block.

\begin{Theorem}\label{thm:2.14}
If there exists an
$(r,\lambda)$-design with $v$ points and $b$ blocks, and a $(b,s \,; \eta)$-DM over $\mathbf{E}_{b}$,
then we have a GDD-type spanning bipartite block design SBBD$^*( s, v, N \,; \Lambda)$, where
$$N=\eta b^2, \quad
\Lambda=(\mu, \lambda_{12}, \lambda_{21}, \lambda_{22})
= ( \eta br, \,
 \eta b \lambda, \,
 \eta r^2, \,  \eta r^2).
$$
It has a $(\eta b^2 \times  sv)$-design matrix. If $s > b-r$, then it satisfies the spanning condition.
\end{Theorem}

\begin{Proof}
Let $H$ be the incidence matrix of an $(r,\lambda)$-design with $v$ points and $b$ blocks.
Let $T_{x_0}, T_{x_1}, \ldots, T_{x_{b-1}},$ $ x_i \in \mathbf{E}_{b}$ be tile matrices  which are defined in Equation (\ref{eq:T}) from $H$.
Suppose $D=[d_{i,j}]$ is an $(\eta b  \times s)$-matrix of $(b,s ; \eta)$-DM over $\mathbf{E}_{b}$.
Let $X$ be the information matrix from  $T_{x_0}, T_{x_1}, \ldots, T_{x_{b-1}}$  using Equations (\ref{def_X}) and (\ref{def_Xj}).
A diagonal submatrix of the information matrix $X^tX$ is from Equation (\ref{pr:TtT}):
$$
X_j^{t} X_j =  \sum_{i=1}^{\eta b} T_{d_{i,j}}^{\,t} T_{d_{i,j}} = \eta b \cdot ( r I_v + \lambda  (J_v-I_v) ) \mbox{\ \ for any  } 1\le i\le s.
$$
Next, consider an off-diagonal submatrix  $X_j^t X_{j'}$, $j \ne j'$.
Let $L_d=\{(x,y)\}$ be a set of pairs such that every difference $d= x-y, d \in \mathbf{E}_b$ occurs exactly $\eta$ times.
From Lemma \ref{lem:diff}, we have
$$
X_j^t X_{j'} = \sum_{(x,y)\in L_d} T_x^t \, T_y = \eta \sum_{x\in \mathbf{E}_b} T_0^t \, T_x\ = \eta \, r^2 J_v, \ j \ne j'.
$$
\qed
\end{Proof}

Each row of $X$ is an SB-block having a form $(\mathbf{x}_1, \mathbf{x}_2,\ldots, \mathbf{x}_s)$, where each $\mathbf{x}_i$ is a row of $H$. If these $\mathbf{x}_i$'s consist of the all  vectors of $H$,  $s=b$,
then each row of $\sum_{i=1}^s X_i$ is $(r,r,\ldots,r)$. For the spanning condition, zeros must not occur in the vector.
If at least $b-(r-1)$ different rows of $H$ appear in an SB-block, then the spanning condition is guaranteed.
When the spanning condition is not satisfied,  we show Example \ref{ex:3} to adjust a difference matrix
using Property \ref{pr:add}.
In the following example, an adjustment  of a difference matrix over $\mathbf{{F}_{2}}^3$ will be seen.

\begin{Example}\label{Ex:BIBDCMD}
Consider a $(4,2)$-design with 7 points and 8 blocks
\[
\{\{1,3,5\}, \{0,3,4\}, \{2,3, 6\}, \{0,1,2\}, \{1,4,6\}, \{0,5,6\}, \{2,4,5\}, \{0,1,2,3,4,5,6\}\}.
\]
Let $\mathbf{E}_b = \mathbf{F}_2 \times \mathbf{F}_2 \times \mathbf{F}_2$.
The incidence matrix is expressed as
\begin{center}
{\scalebox{0.9}[0.9]{
$
 H=
 \begin{bmatrix}
 \mathbf{h}_{(0,0,0)}\\
 \mathbf{h}_{(0,0,1)}\\
 \mathbf{h}_{(0,1,0)}\\
 \mathbf{h}_{(0,1,1)}\\
 \mathbf{h}_{(1,0,0)}\\
 \mathbf{h}_{(1,0,1)}\\
  \mathbf{h}_{(1,1,0)}\\
  \mathbf{h}_{(1,1,1)}\\
 \end{bmatrix}
 =
 \begin{bmatrix}
 & 1 &  & 1 &  & 1 & \\
 1&  &  & 1 & 1 &  & \\
 &  & 1 & 1 &  & & 1\\
 1& 1 & 1 &  &  &  & \\
  & 1 &  &  & 1 &  & 1\\
1 &  &  &  &  & 1 & 1\\
 &  & 1 &  & 1 & 1 & \\
 1 & 1 & 1 & 1 & 1 & 1 & 1
\end{bmatrix}
$.
} }
\end{center}
Then it holds $H^t  H = 4 I_7 +  2 (J_7 - I_7)$.

Using Equation (\ref{eq:T}), the tile matrices
$T_{(0,0,0)}, T_{(1,0,0)},\ldots, T_{(1,1,1)}$ are as follows:
\begin{center}
{\scalebox{0.8}[0.8]{
$
T_{(0,0,0)}=H ,
\quad
 T_{(1,0,0)}=\begin{bmatrix}
  & 1 &  &  & 1 &  & 1\\
1 &  &  &  &  & 1 & 1\\
 &  & 1 &  & 1 & 1 & \\
 1 & 1 & 1 & 1 & 1 & 1 & 1\\
  & 1 &  & 1 &  & 1 & \\
  1&  &  & 1 & 1 &  & \\
  &  & 1 & 1 &  & & 1\\
  1& 1 & 1 &  &  &  &
 \end{bmatrix}
 ,
\quad
T_{(0,1,0)}=\begin{bmatrix}
 &  & 1 & 1 &  & & 1\\
  1& 1 & 1 &  &  &  & \\
 & 1 &  & 1 &  & 1 & \\
 1&  &  & 1 & 1 &  & \\
  &  & 1 &  & 1 & 1 & \\
 1 & 1 & 1 & 1 & 1 & 1 & 1\\
   & 1 &  &  & 1 &  & 1\\
1 &  &  &  &  & 1 & 1
\end{bmatrix},
\ldots.
$
} }
\end{center}

From Property \ref{ppMD}, there is an $(8,8;1)$-DM over $\mathbf{F}_2 \times \mathbf{F}_2 \times \mathbf{F}_2$.
The following difference matrix $D$ is basically from the multiplication table over $\mathbf{F}_{2^3}$, and the $6$-th, $7$-th, and $8$-th columns are added by $(1,0,0), (1,1,0)$, and $(1,0,1)$, respectively.
\begin{center}
{\scalebox{0.8}[0.8]{
$
D =
\begin{bmatrix}
 (0,0,0) & (0,0,0) & (0,0,0) & (0,0,0) & (0,0,0) & (1,0,0)& (1,1,0) & (1,0,1) \\
 (0,0,0) & (1,0,0) & (0,1,0) & (1,1,0) & (0,0,1) & (0,0,1) & (1,0,1) &  (0,1,0)\\
 (0,0,0) & (0,1,0) & (0,0,1)  & (0,1,1) & (1,1,0) & (0,0,0) & (0,0,1) & (0,0,0)\\
 (0,0,0) & (1,1,0) & (0,1,1)  &  (1,0,1) & (1,1,1) & (1,0,1) & (0,1,0) & (1,1,1)\\
 (0,0,0) & (0,0,1) & (1,1,0) &  (1,1,1) & (0,1,1) & (1,1,0) & (0,1,1) & (0,0,1)\\
 (0,0,0) & (1,0,1) & (1,0,0)  & (0,0,1) & (0,1,0) & (0,1,1) & (0,0,0) &  (1,1,0) \\
 (0,0,0) & (0,1,1) & (1,1,1)  &  (1,0,0) & (1,0,1) & (0,1,0) & (1,0,0) & (1,0,0) \\
 (0,0,0) & (1,1,1) & (1,0,1)   &  (0,1,0)  & (1,0,0)  & (1,1,1)  &  (1,1,1) & (0,1,1)
\end{bmatrix}
$
}}
\end{center}
By pasting the tile matrices $T_{(0,0,0)}, T_{(1,0,0)}, \ldots , T_{(1,1,1)}$ into the above difference matrix,
we have a $64 \times 56$ design matrix $X$,  and the following $56 \times 56$ information matrix:
\begin{center}
 {\scalebox{0.8}[0.8]{
$
X^t X =  I_{8} \otimes
\left[\begin{array}{ccccccc}
32 & 16 & 16 & 16 & 16 & 16 & 16  \\
16 & 32 & 16 & 16 & 16 & 16 & 16 \\
16 & 16 & 32 & 16 & 16 & 16 & 16   \\
16 & 16 & 16  & 32 & 16 & 16 & 16   \\
16 & 16 & 16 & 16 & 32 &16 & 16  \\
16 & 16 & 16 & 16 & 16& 32 & 16  \\
16 & 16 & 16 & 16 & 16 & 16 & 32
\end{array}\right]
+(J_{8}-I_{8})  \otimes
\left[\begin{array}{cccccccc}
16 & 16 & 16 & 16 & 16 & 16 & 16 \\
16 & 16 & 16 & 16 & 16 & 16 & 16 \\
16 & 16 & 16 & 16 & 16 & 16 & 16 \\
16 & 16 & 16 & 16 & 16 & 16 & 16 \\
16 & 16 & 16 & 16 & 16 & 16 & 16 \\
16 & 16 & 16 & 16 & 16 & 16 & 16 \\
16 & 16 & 16 & 16 & 16 & 16 & 16
\end{array}\right].
$
}
}
\end{center}
\end{Example}

Table \ref{tbl:primeBIBD} is  a list of existing BIBDs with $b$  blocks, where $b$ is a prime power less than $100$ selected from the table in Colbourn and Dinitz~\cite{Handbook}.
From Property \ref{ppMD}, there exists a $(b,b,1)$-DM over the group $\mathbf{E}_b$. We can construct a GDD type SBBD($b, v, b^2 ; \Lambda), \Lambda=(br, \, b\lambda, \, r^2, \, r^2)$.

\begin{table}[ht]
\begin{center}
\begin{tabular}{ccccc|c}
$v$ & $b$ & $r$ & $k$ & $\lambda$  & \text{Remark} \\ \hline \hline
 7 & 7 & 3 & 3 & 1 & \text{PG(2,2)} \\
 11 & 11 & 5 & 5 & 2 &  \\
 13 & 13 & 4 & 4 & 1 & \text{PG(2,3)} \\
 19 & 19 & 9 & 9 & 4 &  \\
 23 & 23 & 11 & 11 & 5 & \\
 25 & 25 & 9 & 9 & 3 &  \\
 27 & 27 & 13 & 13 & 6 & \text{27=$3^3$} \\
 31 & 31 & 6 & 6 & 1 & \text{PG(2,5)} \\
 31 & 31 & 10 & 10 & 3 &  \\
 31 & 31 & 15 & 15 & 7 & \text{PG(4,2)} \\
 37 & 37 & 9 & 9 & 2 & \\
 41 & 41 & 16 & 16 & 6 &  \\
 43 & 43 & 21 & 21 & 10 & \\
  47 & 47 & 23 & 23 & 11 & \\
 \hline
\end{tabular}
\quad
\quad
\quad
\begin{tabular}{ccccc|c}
$v$ & $b$ & $r$ & $k$ & $\lambda$  & \text{Remark} \\ \hline \hline
   7 & 49 & 21 & 3 & 7 &  \\
 49 & 49 & 16 & 16 & 5 & \text{49=$7^2$} \\
 59 & 59 & 29 & 29 & 14 &  \\
 61 & 61 & 16 & 16 & 4 & \\
 61 & 61 & 25 & 25 & 10 &  \\
 67 & 67 & 33 & 33 & 16 &  \\
 71 & 71 & 15 & 15 & 3 & \\
 71 & 71 & 21 & 21 & 6 & \\
 71 & 71 & 35 & 35 & 17 & \\
 73 & 73 & 9 & 9 & 1 & \text{PG(2,8)} \\
 79 & 79 & 13 & 13 & 2 &  \\
 79 & 79 & 27 & 27 & 9 & \\
 79 & 79 & 39 & 39 & 19 & \\
  & & & & &\\ \hline
\end{tabular}
\end{center}
  \caption{BIBD with prime power $b$ blocks}
  \label{tbl:primeBIBD}
\end{table}

\section{Decomposition method}

Let ($V, \mathcal{B}$) be an  $(r,\lambda)$-design with $v$ points and $b$ blocks.
When $b>>v$, the method described in Section 3 can only construct SBBDs in which $v_1$ and $v_2$ are significantly different.
In this section, we propose a construction method to meet the requirement to have SBBDs in which $v_1$ and $v_2$ are relatively close using an $(r,\lambda)$-design with $b>>v$.

Let $ \mathcal{B}_1, \mathcal{B}_2, \ldots, \mathcal{B}_m $ be a
partition of  $\mathcal{B}$,  where $ |\mathcal{B}_i | = b_i $ and each point of $V$ is contained in $\mathcal{B}_i$ at least one for $i=1,2,\ldots,m$.
Let $\mathbf{E}_{b_i}^{(i)}$ be a group of order $b_i$, $1\le i \le m$.
Then the $(b_i \times v)$-incidence matrix $H_i$ between  $\mathcal{B}_i$ and $V$ is described as
$$
H_i =
\begin{bmatrix}
\mathbf{h}_{x_0}^{(i)}\\
\mathbf{h}_{x_1}^{(i)}\\
\vdots\\
\mathbf{h}_{x_{b_i -1}}^{(i)}
\end{bmatrix}
, \mbox{ where } x_j\in \mathbf{E}_{b_i}^{(i)} \mbox{ for } 1 \leq i \leq m.
$$

\begin{Property}\label{HtH}
If $(V,\mathcal{B})$ is an $(r,\lambda)$-design with $v$ points and $b$ blocks,
 then
 $$\sum_{i=1}^m H_i^t H_i = r I_v + \lambda(J_v-I_v).$$
\end{Property}
For each $H_i$, $1 \leq i \leq m$, we generate $b_i$ tile matrices $T_y^{(i)}, y \in \mathbf{E}_{b_i}^{(i)}$,
of the size  $(b_i \times v)$ by adding an element of $ \mathbf{E}_{b_i}^{(i)}$ to the subscripts, same as Equation (\ref{eq:T}).
Let
\begin{equation}\label{eq: T2}
T_y^{(i)}=
 \begin{bmatrix}
\mathbf{h}_{x_0 +y }^{(i)}\\
\mathbf{h}_{x_1 +y}^{(i)}\\
\vdots\\
\mathbf{h}_{x_{b_i -1} +y}^{(i)}
\end{bmatrix}
, \mbox{ where } y\in \mathbf{E}_{b_i}^{(i)}
\end{equation}

It is not difficult to see the following for the tile matrices $T_y^{(i)}$, $y \in \mathbf{E}_{b_i}^{(i)}$:
\begin{itemize}
\item $T_y^{(i)}$ is a $b_i \times v$ matrix for any  $y \in \mathbf{E}_{b_i}^{(i)}$,
\item  $(T_y^{(i)})^t \, T_y^{(i)}  = H_i^t H_i$ for any  $y \in \mathbf{E}_{b_i}^{(i)}$.
\item  the set of rows of  $T_y^{(i)}$  is exactly equal to the set of rows of $H_i$ for any $y \in \mathbf{E}_{b_i}^{(i)}$,
which implies
\begin{gather}
  \sum_{i=1}^m (T_{y}^{(i)})^t \, T_{y}^{(i)} = \sum_{i=1}^m {H_i}^t H_i =r I_v + \lambda (J_v-I_v) \   \mbox{ for any } y \in \mathbf{E}_{b_i}^{(i)}. \label{eq:tile7}
  \end{gather}
\end{itemize}

\begin{Lemma}\label{lem:tilepro3}
Let $(V, \mathcal{B})$ be an $(r,\lambda)$-design and
let  $T_y^{(i)}$,  $1\le i \le m$, $y \in \mathbf{E}_{b_i}^{(i)}$, be the tile matrices of $\mathcal{B}_i$, $|\mathcal{B}_i|=b_i$,
where $\{\mathcal{B}_1, \mathcal{B}_2,\ldots, \mathcal{B}_m \}$ is a partition of $\mathcal{B}$.
If  every element of $V$ appears in $\mathcal{B}_i$ exactly $r_i$ times for $1\le i\le m$,
then
the following two equations hold:
\begin{align}
& (T_x^{(i)})^t \,T_{y}^{(i)} =   (T_{x+d}^{(i)})^t \, T_{y+d}^{(i)} \
  \mbox{ for any }  x, y, d \in \mathbf{E}_{b_i}^{(i)}, \, 1\le i\le m, \label{eq:tile8} \\
& \sum_{z \in \mathbf{E}_{b_i}^{(i)}} (T_y^{(i)})^t  \, T^{(i)}_{y+z} =  {r_i}^2 J_v \  \mbox{ for any } y \in \mathbf{E}_{b_i}^{(i)}, \, 1\le i\le m,\label{eq:tile9}
\end{align}
\end{Lemma}

\begin{Proof}
First, we prove Equation (\ref{eq:tile8}) from Equation (\ref{eq:diff}).
Suppose that $H$ is divided into $H_1, H_2,\ldots$, $H_m$, and  every column of $H_i$ has $r_i$ ones.
Let  $W_i = \{T_{x+z}^{(i)}\, | \, z\in  \mathbf{E}_{b_i}^{(i)}\}$ be the set of tile matrices produced from $H_i$.
Then each row of $H_i$ appears exactly once in the same rows of tile matrices in $W_i$.
Using the same approach as in the proof of (\ref{eq:diff}), we  have
$$
(T_x^{(i)})^t \,T_{y}^{(i)} =   (T_{x+d}^{(i)})^t \, T_{y+d}^{(i)} \
  \mbox{ for any }  x, y, d \in \mathbf{E}_{b_i}^{(i)}, \, 1\le i\le m.
  $$
Next, Equation (\ref{eq:tile9}) can be proved in the same manner as the proof of Equation (\ref{eq:diff3}),
$$
\sum_{z \in \mathbf{E}_{b_i}^{(i)}}   (T_y^{(i)})^t \, T^{(i)}_{y+z}
= {r_i}^2 J_v  \mbox{\ \  for } y \in \mathbf{E}_{b_i}^{(i)}.
$$
\qed
\end{Proof}

Let $D^{(i)}=\left[d_{pq}^{(i)}\right]$ be
$(b_i, s ; \eta_i)$-DM over $\mathbf{E}_{b_i}^{(i)}$ of
size $(\eta b_i  \times s)$, $1 \leq i \leq m$.
We paste the tile matrices $T_1^{(i)}, T_2^{(i)}, \ldots, T_{b_i}^{(i)}$
on the difference matrix $D^{(i)}=\left[d_{pq}^{(i)}\right]$, and denote it by
\begin{gather}
 X^{(i)} = \left[\, T_{d_{pq}^{(i)}}^{(i)} \,\right]. \label{eq:tile11}
\end{gather}
Then we have an $(\eta \sum_{i=1}^m b_i^2 \times s v)$-design matrix
\begin{gather}
 X=
\begin{bmatrix}X^{(1)}\\ X^{(2)}\\ \vdots \\  X^{(m)}\end{bmatrix}
=( X_1 | X_2| \cdots | X_s ).  \label{eq:tile12}
\end{gather}

\begin{Theorem}
If there is an $(r,\lambda )$-design $(V, \mathcal{B})$  with $v$ points and $b$ blocks
which is partitionable into $\mathcal{B}_1, \mathcal{B}_2,\ldots, \mathcal{B}_m$ such that every point of $V$ appears in $\mathcal{B}_i$ exactly $r_i$ times,
and if  there exist  $(b_i,s; \eta)$-difference matrices, $i=1,2,\ldots, m$, satisfying
$b=b_1+b_2+\cdots+b_m$,
then
there exists a GDD-type SBBD $(s, v, N ; \Lambda)$, where $N= \eta \sum_{i=1}^m b_i^2$ and
$$
\Lambda = (\mu, \lambda_{12}, \lambda_{21}, \lambda_{22})
 =(\eta r b, \, \eta \lambda b, \, \eta \sum_{i=1}^m {r_i}^2, \, \eta  \sum_{i=1}^m {r_i}^2).
 $$
\end{Theorem}

\begin{Proof}
First,  we compute the diagonal submatrix $X_j^t X_j$ of $X^tX$.
From  Equation (\ref{eq:tile7}), we have
$$
X_j^t \, X_j =  \sum_{i=1}^{m} \sum_{p=1}^{\eta b_i}T_{d_{p,j}^{(i)}}^{(i) \, t} \, T_{d_{p,j}^{(i)}}^{(i)} = \eta b \cdot ( r I_v + \lambda  (J_v-I_v) ) \mbox{  for any  } 1\le i\le s.
$$
Second, we compute an off-diagonal submatrix $X_j^t X_{j'}, \ 1 \le  j\ne j' \le s.$
The following equation holds regardless of the elements of $x_i \in \mathbf{E}_{b_i}^{(i)}$, $1\le i \le m$, from Equations (\ref{eq:tile8}) and (\ref{eq:tile9}):
$$
X_j^t X_{j'} =   \eta \sum_{i=1}^m \sum_{z \in \mathbf{E}_{b_i}^{(i)} }  (T_{x_i}^{(i)})^t T_{x_i+z}^{(i)}
 = \eta \, \sum_{i=1}^m r_i^2 J_v.
 $$
\qed
\end{Proof}

Suppose we want to have SBBDs of $K_{v_1,v_2}$ such that  $v_1$ and  $v_2$ are as close as possible.
Let $(V, \mathcal{B})$ be an $(r,\lambda)$-design with $v$ points and $b$ blocks,
and let $\mathcal{B}_1, \mathcal{B}_2,\ldots, \mathcal{B}_m$ be a partition of $\mathcal{B}$.
 When decomposing the block set, the following should be considered:
\begin{itemize}
    \item   $b=b_1+b_2,+\cdots +b_m$, where $b_i=|\mathcal{B}_i|$, $1 \leq i \leq m$,
    \item  every point of $V$ appears in $\mathcal{B}_i$ exactly $r_i$ times,  $1 \leq i \leq m$,
    \item  each $b_i$ is as close to $v (= v_2)$ as possible,
    \item  each $b_i$ is a prime or prime power (When it is hard to decompose into  such $b_i$'s,  we can have new $(r+1,\lambda+1)$-design with $b+1$ blocks by adding a block $B_{b+1} = V$),
    \item  an integer $s (=v_1)$ in $(b_i,s;\eta)$-DM is
          $s= \min\{b_1,b_2,\ldots, b_m\} $.
\end{itemize}

\begin{Example}
Consider a $(5,3,3)$-BIBD with $10$ blocks. The set of  blocks is divided into two parts, each  consisting of 5 blocks.
Their incidence matrices $H_1, H_2$ of those  two parts are as follows:

\begin{center}
$ H_1 =
\begin{bmatrix}
  0 &  0  &  1  &  1  &  1 \\
  1 &   0  &  0  &  1  &  1 \\
  0 &  1  &  1   & 1  &  0 \\
   1 &   1 &   0 &   0   & 1\\
  1  &  1  &  1  & 0  &  0
\end{bmatrix}
, \ H_2=
\begin{bmatrix}
  1 &  1  &  0  &  1  &  0 \\
  1 &   0  &  1  &  1  &  0 \\
  0 &  1  &  1   & 0  &  1 \\
   1 &   0 &   1 &   0   & 1\\
  0  &  1  &  0  & 1  & 1
\end{bmatrix}
$.
 \end{center}
Naturally,  ${H_1}^t H_1 + {H_2}^t H_2= 6 I_5 + 3 (J_5-I_5)$.
 Since $b_1=b_2=5$, there exist the following difference matrices $D^{(1)}$, $D^{(2)}$ over the group  $\mathbf{Z}_5 =\{0,1,2,3,4\}$:
 $$
 D^{(1)}=
 \begin{bmatrix}
   0  &  0  &  1  &  4  &  3\\
   0  & 1  &  3  &  2  &  2\\
   0  &  2 &   0 &   0 &   1\\
   0  &  3 &   2 &  3 &   0\\
   0  &  4 &   4  &  1  &  4
 \end{bmatrix}
  , \ D^{(2)}=
 \begin{bmatrix}
   0  &  0  &  1  &  4  &  3\\
   0  & 1  &  3  &  2  &  2\\
   0  &  2 &   0 &   0 &   1\\
   0  &  3 &   2 &  3 &   0\\
   0  &  4 &   4  &  1  &  4
 \end{bmatrix}
 .$$
 From $H_1$, we can produce tile matrices $T_0^{(1)}, T_1^{(1)}, \ldots, T_4^{(1)}$ by the method of Equation (\ref{eq: T2}),
 \begin{center}
\scalebox{0.9}[0.9]{
$ T_0^{(1)} = H_1
, \ T_1^{(1)}=
\begin{bmatrix}
  1 &   0  &  0  &  1  &  1 \\
  0 &  1  &  1   & 1  &  0 \\
   1 &   1 &   0 &   0   & 1\\
  1  &  1  &  1  & 0  &  0\\
    0 &  0  &  1  &  1  &  1
\end{bmatrix}
, \ T_2^{(1)}=
\begin{bmatrix}
  0 &  1  &  1   & 1  &  0 \\
   1 &   1 &   0 &   0   & 1\\
  1  &  1  &  1  & 0  &  0\\
    0 &  0  &  1  &  1  &  1 \\
  1 &   0  &  0  &  1  &  1
\end{bmatrix}
,\ldots,
\ T_4^{(1)}=
\begin{bmatrix}
  1  &  1  &  1  & 0  &  0 \\
  0 &  0  &  1  &  1  &  1 \\
  1 &   0  &  0  &  1  &  1 \\
  0 &  1  &  1   & 1  &  0 \\
   1 &   1 &   0 &   0   & 1
\end{bmatrix}
, $}
\end{center}
and similarly $T_0^{(2)}, T_1^{(2)}, \ldots, T_4^{(2)}$ from $H_2$.
Finally, we paste the tile matrices $T_i^{(1)}$ and $T_i^{(2)}$ onto the difference matrices $D^{(1)}$ and $D^{(2)}$, respectively.
Then we have a GDD type SBBD($5,5, 50; \Lambda)$, where $\Lambda=(30,15,18,18)$.
 Its information matrix is
 $$
 X^t X = I_5 \otimes
 \begin{bmatrix}
 30 & 15 & 15 & 15 & 15 \\
 15 & 30 & 15 & 15 & 15 \\
 15 & 15 & 30 & 15 & 15 \\
 15 & 15 & 15 & 30 & 15 \\
 15 & 15 & 15 & 15 & 30
 \end{bmatrix}
 + (J_5-I_5) \otimes
  \begin{bmatrix}
18 & 18 & 18 & 18 & 18 \\
18 & 18 & 18 & 18 & 18 \\
18 & 18 & 18 & 18 & 18 \\
18 & 18 & 18 & 18 & 18 \\
18 & 18 & 18 & 18 & 18
  \end{bmatrix}
\  .$$
\end{Example}

\section{Optimal design and existence}

Takeuchi~\cite{Takeuchi1961} shows that a  specific type group divisible design is optimum in a statistical model.
In this section, we discuss statistical models of  the GDD-type SBBD and optimality.

Let $\mathbf{y}$, $\boldsymbol{\tau}$, and $\boldsymbol{\epsilon}$ be vectors of data, main effects, and errors, respectively.  $\mu$ is the central effect.
$X = [x_{ij} ]$ is an $N \times v$  $(0, 1)$-matrix called a design matrix.
Each data is obtained as the sum of some effects.
Then the model can be represented as
\begin{equation}\label{model}
\begin{split}
&\mathbf{y} = \mu\mathbf{1}_N +  X \boldsymbol{\tau} + \boldsymbol{\epsilon} , \\
&\boldsymbol{\tau}^t\mathbf{1}_v =0 .
\end{split}
\end{equation}

When evaluating efficient designs, the smaller the variance of the estimator, the better.
The goodness is quite different for design matrices of the same size and the same number of 1s.
Since there is usually  more than one estimator,
there are several criteria for design optimality.
Here we show a  criterion of optimality called {\it E-optimal}.

\begin{Definition}[E-optimality, Kiefer~\cite{Kiefer1958}, Shar and Sinha~\cite{Shar-Sinha1960}]
Let $\Omega$ be a class of $N\times v$ (0,1)-matrices $X$ having the same number of ones. If the following function has a maximum value for $X$ in $\Omega$, then the design matrix $X$ is called  E-optimum relative to $\Omega$:
$$
\min_{1\le i \le v-1} \{\theta_i\},
$$
where $\theta_1, \theta_2,\ldots, \theta_{v-1}$ , $\theta_i > 0$, are the eigenvalues of $X^{t}X$.
\end{Definition}

The optimality of group divisible designs is discussed in Takenchi~\cite{Takeuchi1961}. 
The statistical models for group divisible designs do not consider the group structure of the variety set $V$ at all. That is the model of (\ref{model}) is assumed with $N$ blocks and $v = mg$.
Let $\Omega$ is the class of $N\times v$ $(0,1)$-matrices $X$ which contain exactly  $kN$ ones.
\begin{Theorem}[Takeuchi~\cite{Takeuchi1961, Takeuchi1963}]
A group divisible design GD$_{\lambda_1,\lambda_2} (k,  g \,;\, v)$ with $\lambda_2=\lambda_1 + 1$ is E-optimum relative to $\Omega$.
\end{Theorem}
Naturally, this theorem applies to SBBDs constructed from the GDDs.

\begin{Theorem}
    An GDD type SBBD$^*$($v_1,v_2, N ; \Lambda)$, $\Lambda=(r, \lambda_1, \lambda_2, \lambda_2 )$, where $v=v_1v_2$ and $\lambda_2=\lambda_1 +1 $ is E-optimum relative to $\Omega$.
\end{Theorem}

Many group divisible designs with $\lambda_2 =\lambda_1 + 1$ are known. We introduce some well-known constructions in this section.
Suppose here  $(V,\mathcal{B})$ be a  ($v,k,1$)-BIBD.
 Let $\Pi = \{ \Pi_1, \Pi_2, \ldots, \Pi_n\}$, $|\Pi_i|=g$, be  a partition of $V$.
 Each  block of $\mathcal{B}$ intersects $\Pi_i$ at no points, one point, or at all points in the block.
 If $\mathcal{B}' \subset \mathcal{B}$ consists of the blocks not contained in any $\Pi_i$, then $(V,\, \mathcal{B'})$ is a  $GD_{0,1}(k, g \,; v)$.

\begin{Example}
The points and the lines of PG$(n,q)$,  $q$ a prime power,  form a $( (q^{n+1} - 1 )/ (q-1),q+1,1)$-BIBD.
There exists a parallel class of $t$-flats (equivalently a $t$-spread of PG($t,q$)) if and only if  $(t+1) \mid (n+1)$.
That is, there exists a  $GD_{0,1}(q+1, (q^{t+1}-1)/(q-1) ; v)$.
\end{Example}

\begin{Example}
    In AG($n,q$), $q$ a prime power, there is a parallel class of t-flats $\Pi_i$ (= AG($t,q$)) for $1\le t \le n-1$. The points and the lines is a ($q^n, q, 1$)-BIBD 
    with parallel class of $\Pi_i,\  |\Pi_i|=g=q^t$. That is, there is a $GD_{0,1}(q, q^t ; q^n)$.
\end{Example}

\begin{Example}
    For any $q$ prime power, there is an orthogonal array GD$_{0,1}( q+1, q  ; q(q+1))$.
\end{Example}

\begin{Definition}[Complement design]
The complement design of $(V,\mathcal{B})$ is $(V,\overline{\mathcal{B}})$, where $\overline{\mathcal{B}}= \{ V\setminus B \  |\, B\in \mathcal{B} \}$.
\end{Definition}
\begin{Property}
The complement design of $GD_{\lambda_1,\lambda_2}(k, v_2;v_1v_2)$ with $N$ blocks is $GD_{\lambda_1', \lambda_2'} (v-k, v_2;$ $v_1v_2)$, where
$\lambda_1' = N-2r+\lambda_1$, $\lambda_2'=N-2r+\lambda_2$ and $r=kN / (v_1v_2)$.
Therefore if $\lambda_2=\lambda_1 + 1$ then $\lambda_2'=\lambda_1' + 1$.
 \end{Property}

\begin{Property}
    Let $V$ be the point set of $AG(n,q)$, $n \ge 2$,  $q$ a prime power. The set of groups $G_1, G_2,\ldots, G_q$  is a parallel class of hyperplanes. The blocks $\mathcal{B}$ is the set of hyperplanes not any of $G_i$.
    Then $(V, \mathcal{B})$ is a $GD_{\lambda_1,\lambda_2}(q^{n-1}, q^{n-1} ;\, q^n)$, where
    $$\lambda_1 = \frac{q^{n-1}-q}{q-1}  , \ \  \lambda_2 =\frac{q^{n-1}-1}{q-1}.$$
    That is, $\lambda_2 = \lambda_1 +1 $.
\end{Property}
 \Proof
 Any hyperplane of $\mathcal{B}$ meets each $G_i$ in a $(n-2)$-flat and these $(n-2)$-flats are parallel.
 Suppose $p_1, p_2$ are two points in $G_1$.
 There are $(q^{n-2}-1)/(q-1)$ $(n-2)$-flats containing $p_1$ and $p_2$.
 Let $\pi$ be one of them.
 Consider hyperplanes  of $\mathcal{B}$ containing $\pi$.
 There are $q$ $(n-2)$-flats in $G_2$ parallel to $\pi$.
 One flat of them and $\pi$ determine a unique hyperplane of $\mathcal{B}$.
 So, there are $\lambda_1=(q^{n-2}-1)/(q-1)\times q$ hyperplanes of $\mathcal{B}$ containing $p_1$ and $p_2$.
 Next, consider two points $p_1$ in $G_1$ and $p_2$ in $G_2$.
 There are $(q^{n-1}-1)/(q-1)$ $(n-2)$-flats in $G_1$ containing $p_1$.
 An $(n-2)$-flat in $G_1$ and the point $p_2$ in $G_2$ determines a unique hyperplane of $\mathcal{B}$.
 That is, $\lambda_2= (q^{n-1}-1)/(q-1)$.
 \qed

\subsection*{Acknowledgments}
This work was supported by JSPS KAKENHI Grant Numbers JP19K11866 and 21K13845.

\bibliographystyle{abbrv}
\bibliography{ref}

\end{document}